\documentclass[10pt]{amsart}
\usepackage[english]{babel}
\usepackage{amsmath}
\usepackage{amssymb}
\usepackage{amsfonts}
\usepackage{graphicx}

\newtheorem{lemma}{Lemma}
\newtheorem{theorem}{Theorem}

\newtheorem{corollary}{Corollary}
\newtheorem{proposition}{Proposition}

\author{Alexei Ya. Kanel-Belov}
\address{College of Mathematics and Statistics, Shenzhen University, Shenzhen, 518061, China}
\email{kanelster@gmail.com}

\author{Igor Melnikov}
\address{Moscow Institute of Physics and Technology, Dolguprudny, Russia}
\email{melnikov\_ig@mail.ru}

\author{Ivan Mitrofanov}
\address{C.N.R.S., \'{E}cole Normale Superieur, PSL Research University, France }
\email{phortim@yandex.ru}

\title{On cogrowth function of algebras and its logarithmical gap\footnote{The paper was supported  by Russian Science Foundation (grant no. 17-11-01377)}}

\begin{document}

\maketitle

\begin{abstract}

Let  $A \cong k\langle X \rangle / I$ be an associative algebra.
A finite word over alphabet $X$ is $I${\it-reducible} if its image in $A$ is a $k$-linear combination of length-lexicographically lesser words.
An {\it obstruction} in a subword-minimal $I$-reducible word.
If the number of obstructions is finite then $I$ has a finite Gr{\"o}bner basis, and the word problem for the algebra is decidable.
A {\em cogrowth} function is number of obstructions of length $\le n$.
We show that the cogrowth function of a finitely presented algebra is either bounded or at least logarithmical.
We also show that an uniformly recurrent word has at least logarithmical cogrowth.
	
\end{abstract}

\begin{abstract}
Soit $A\cong k\langle X\rangle/I$ une alg{\`e}bre associative.
Un mot fini sur l'alphabet $ X $ est $ I $ {\ it-reductible} si son image dans $ A $ est une combinaison lin{\'e}aire $ k $ de mots de longueur lexicographiquement moindre.
Une {\it obstruction} dans un mot minimal $I$ -r{\'e}ductible.
Si le nombre d'obstructions est fini, alors $I$ a une base finie Gr{\"o}bner, et le mot probl{\`e}me pour l'alg{\`e}bre est d{\'e}cidable.
Une fonction {\em co-croissance} est le nombre d'obstructions de longueur $\le n $.
Nous montrons que la fonction de co-croissance d'une alg{\`e}bre finement pr{\'e}sent{\'e}e est soit born{\'e}e, soit au moins logarithmique.
Nous montrons {\'e}galement qu'un mot uniform{\'e}ment r{\'e}current a au moins une co-croissance logarithmique.
\end{abstract}

\section{Cogrowth of associative algebras}.

Let $A$ be a finitely generated algebra over a field $k$.
Then $A \cong k\langle X \rangle / I$, where $k\langle X \rangle$ is a free algebra with generating set $X = \{x_1,\dots, x_s\}$  and $I$ is a two-sided {\it ideal of relations}.
Further we assume the generating set is fixed.
Let ``$\prec$'' be a well-ordering of $X$,
$x_1 \prec \dots \prec x_s$.
This order can be extended to a linear order on the
set $\langle X \rangle$ of monomials of $k\langle X \rangle$: $u_1 \prec u_2$ if $|u_1| < |u_2|$ or $|u_1| = |u_2|$ and $u_1 <_{lex} u_2$.
Here $|\cdot|$ denotes the length of a word, i.e. the degree of a monomial, and $<_{lex}$ is the lexicographical order.
We denote the set of monomials of degree at most $n$ by $\langle X \rangle_{\leq n}$.
For $f \in k\langle X \rangle$ then we denote the leading (with respect to $\prec$) monomial of $f$ by $\hat{f}$.
If $I$ is a finitely generated ideal, the algebra $k\langle X \rangle / I$ is called {\it finitely presented}.

The {\it growth} $V_A(n)$ is the dimension  $\dim(\text{span}(A_n))$, where $A_n$ is the set of images of $\langle X \rangle_{\leq n}$ in $A$.
We call a monomial $w\in \langle X \rangle$ {\it $I$-reducible} if $w = \hat{f}$ for some relation $f\in I$.
It is easy to see that $V_A(n)$ is equal to the number of $I$-irreducible monomials in $\langle X \rangle_{\leq n}$.

We call a monomial $w\in \langle X \rangle$ an {\it obstruction} in $A$ if $w$ is $I$-reducible, but any proper subword of $w$ is $I$-irreducible.
The {\it cogrowth} of algebra $A$ is defined as the function $O_A(n)$, the number of obstructions of length $\leqslant n$.

A {\it Gr{\"o}ebner basis} of an ideal $I$ is a subset $G \subseteq I$ such that for any $f \in I$ there exists $g \in G$ such that the leading monomial of $f$ contains the leading monomial of $g$ as a subword.

The {\it word problem} for a finitely presented
$k\langle X \rangle$, i.e. the question whether a given element $f\in k\langle X \rangle$ belongs lies in $I$, is undecidable in general case \cite{Ceit, Scott}.
But if $I$ has a finite Gr{\"o}ebner basis $G$, then  $A$ has a decidable {\it word problem}.
If $\hat{f}$ contains $\hat{g}$ for some $g\in G$, then $f$ can be replaced by $f'$ such that $f' - f \in I$ and $\hat{f'} \prec \hat{f}$. This operation is called a {\it reduction}.
After some number of reductions we obtain either 0 or an element $f''$ such that $\hat{f''}$ is $I$-irreducible. In this case, $f \not\in I$.

Note that the problem whether a given element in a finitely presented associative algebra is zero divisor (or is it nilpotent) is undecidable,
even if we are given a finite Gr{\"o}ebner basis
\cite{IvanMal}.

\begin{theorem}
	
	Let $A$ be
	a finitely
	presented algebra and
	let $m$ be the maximum length of its defining relation, $N\ge m$.
	Suppose there are no obstructions of length from the segment
	$[N, 2N]$.
	Then $A$ has a finite Gr{\"o}bner basis.
\end{theorem}

{\it Sketch of proof.}
Let $S$ be the set of all obstructions in $\langle X \rangle_{\leq N}$.
Take for each monomial $w\in S$ a reduced relation $f_w$ such that $\hat{f_w} = w$.
If for some $u_1, u_2, u_3 \in \langle X \rangle$ and
$w_1, {w_2} \in S$ it holds
$u_1u_2 = w_1$ and $u_2u_3 = w_2$, then the word $u_1u_2u_3$ is called a {\it composition} of $f_{w_1}$ and $f_{w_2}$, and the normed element
\[
(f_{w_1} - w_1)u_3 - u_1(f_{w_2} - w_2)
\]

is the {\it result of the composition}.

If we take any two elements of form $f_{w}|w\in S$, the leading monomial of any their composition has length less then $2N$, and can be reduced to zero.
From Bergman's {\it diamond lemma} \cite{Berg} it follows that the set
$\{f_{w_i}|w_i\in S\}$ forms a Gr{\"o}ebner basis for $I$.
$\square$

\begin{corollary}
	Let $A$ be a finitely presented algebra.
	Then the cogrowth function $O_A(n)$ is either constant or no less than logarithmic:
	
	$$O_A(n)\ge \log_2(n)-C.$$
	
	The constant $C$ depends only on the maximal length of a relation.
\end{corollary}

Well known Bergman gap theorem says that the growth function $V_A(n)$  is either constant, linear of no less than
$n(n+3)/2$.

\section{Colength of a period}.

A {\it monomial algebra} is an finitely generated associative algebra whose relations are monomials.
The irreducible monomials of a monomial algebra is the set of all finite words that avoid ``forbidden'' subwords from the list of relations.

For monomial algebras with a finite Gr{\"o}ebner basis, as well as for automaton monomial algebras, the nilpotency problem is algorithmically decidable \cite{BBL, Iyu, KB}, unlike the situation in general case \cite{IvanMal}.

Let $X = \{x_1,\dots, x_k\}$ be a finite alphabet.
We consider infinite sequences on $X$, i.e. maps $X^{\mathbb{N}}$.
A sequence $(a_i)$ is {\it periodic} with period $u$ if
$(a_i) = uuu\dots$.

A finite word $v$ is an {\it obstruction} for a sequence $W$ if $v$ is not a subword of $W$ but any proper subword $v'$ of $v$ is a subword of $W$.
Let $u$ be a finite word. The number of obstructions for $u^{\infty}$ is always finite, we call this number the {\it colength} of the period $u$.
We say that the period is {\it defined by} the set of obstructions.

In \cite{Cheln1}, G.~R.~Chelnokov proved that a sequence of minimal period $n$ can not be defined by less than $\log_2n + 1$ obstructions.
G.~R.~Chelnokov also gave for infinitely many $n_i$ an example of a binary sequence with minimal period $n_i$ and colength of the period $\log_{\varphi}n_i$, where $\varphi = \frac{\sqrt{5}+1}{2}$.

P.~A.~Lavrov found the precise lower estimation for colength of period.

\begin{theorem}\cite{Lavr1}
	Let $A = \{a, b\}$ be a binary alphabet.
	Let $u$ be a word of length $n$ and colength $c$, then $\varphi_c\geq n$, where $\varphi_c$ is the $c$-th Fibonacci number ($\varphi_1=1$, $\varphi_2=2$, $\varphi_3=3$, $\varphi_4=5$ etc.).
\end{theorem}

The case of arbitrary alphabet was considered in \cite{Lavr2} by P.~A.~Lavrov and later in \cite{BogdCheln} by I.~I.~Bogdanov and G.~R.~Chelnokov.

\section {Cogrowth function for an uniformly recurrent sequence}.

Let $A = k\langle X \rangle/F$ be an infinite-dimensional monomial algebra such that adding any new monomial relation to $F$ gives a finite-dimensional algebra.
Irreducible monomials of $A$ are all finite subwords of some {\it uniformly recurrent} sequence \cite{BBL}

A sequence of letters $W$ on a finite alphabet is called
{\it uniformly recurrent} if for any finite subword $u$ of $W$ there exists a number $C(u, W)$ such that any subword of $W$ with length $C(u, W)$ contains $u$.

Again, a finite word $u$ is  an {\it obstruction} for $W$ if it is not a subword of $W$ but any its proper subword is a subword of $W$.
The {\it cogrowth function} $O_W(n)$ is the number of   obstructions with length $\leqslant n$.
The linearly equivalence class of the cogrowth function is an important topological invariant of the corresponding symbolic dynamical system \cite{Beal}.

We prove
\begin{theorem}\label{thm:ur}
	Let $W$ be an uniformly recurrent non-periodic sequence on a binary alphabet.
	Then $$\overline{\lim_{n_i \to \infty}} \frac{O_W(n_i)}{\log_3n} \geq 1.$$
	
\end{theorem}

{\it Sketch of proof.}

The {\it factor language} of $W$ is the set of all its finite subwords.
If $W$ is uniformly recurrent, then the factor language of $W$ is inclusion-minimal among infinite factor languages, i.e. it is not possible to forbid any new word without forbidding all but finite number of words in the language.

This language can be described in terms of {\it Rauzy graphs}. The vertices of the directed graph $R_n(W)$ are subwords of $W$ of length $n$, the edges of $R_n(W)$ are subwords of length $n+1$.
The sequence $(R_n)$ can be constructed by sequential application of operations of two types:

\begin{enumerate}
	\item deleting an edge $H \to H-e$;
	\item $H \to L(H)$, where $L(H)$ is the directed line graph of $H$, (vertices of $L(H)$ are edges of $H$, edges of $L(H)$ are 2-paths in $H$).
\end{enumerate}

For a directed graph $H$ we define its {\it entropy regulator}: $er(H)$ is the minimal integer such that any directed path in $H$ contains at least one vertex with outgoing degree 2.
We show by induction on $k$ that $er(R_n(w)) \leq 2^{O_W(n)}$.

\begin{lemma}\label{lm:del_edge}
	Let $H_0$ be a strongly connected graph, let $H_1$ be $L^{3er(H)}(H)$ and let $e$ be an arbitrary edge in $H_1$.
	Then the digraph $H_1 - e$ contains a strongly connected subgraph $H_2$ such that $er(H_2)\leq 3 er(H)$.
\end{lemma}

Now suppose $O_W(n) < \log_3n$ for all $n > n_0$.
Then we can choose $n_1$ and choose for each obstruction $u_i$ of length $|u_i|>n$ its proper subword $v_i$ such that the sequence of lengths of these subwords is $3K, 9K, \dots, 3^kK$, where $K = er(R_{n_1}(W))$.

Using Lemma \ref{lm:del_edge} we show by induction that for any $n$ the Rauzy graph $R_n(W)$ contains a non empty subgraph $R'_n$ such that vertices of $R'_n$ do not have any of $v_i$ as subwords.

But this contradicts the inclusion  minimality of the factor language of $W$. $\square$

\medskip

Consider a finite alphabet $\{a,b\}$ and the sequence of words $u_i$, defined recursively as
\[u_0 = b, \:u_1 = a, \:u_{k} = u_ku_{k-1} \text{\:for $k \geq 2$}.
\]
the sequence $(u_i)$ has a limit, called {\it Fibonacci word}.
\[F = abaababaabaab\dots
\]

It can be shown that the obstructions of the Fibonacci word $bb, aaa, babab,\dots$
have lengths equal to Fibonacci numbers, so
$O_F(n) \sim \log_{\varphi}n$, where $\varphi = {\frac{\sqrt5 + 1}2}$.

The next propositions shows that in Theorem \ref{thm:ur} we can not replace $\overline{\lim}$ by $\lim$.

\begin{proposition}
There exists an uniformly recurrent non periodic sequence $W$ such that  $\underline{\lim_{n\to \infty}} \frac{O_W(n)}{\ln n} = 0$.
\end{proposition}

{\it Sketch of proof.}
We call a factor language $\mathcal{L}$ {\it uniformly recurrent at level $t$} if for some $T$ for any pair of words $u, U\in \mathcal{L}, |u| = t, |U| = T$ the word $u$ is a subword of $U$.
We can construct a factor language by adding obstructions one by one.
We can wait as long as we want without adding any obstructions to make $\frac{O_W(n)}{\ln n}$  arbitrary small.
After that we can forbid long words to make the factor language uniformly recurrent at some new level, and we iterate these operations infinitely many times.
$\square$

It is easy to see that the cogrowth function of such a word cannot be equal to a cogrowth function of any finitely presented algebra.

\end{document}